\documentclass[10pt]{amsart}

\usepackage{graphicx}
\usepackage{amsmath,amsthm,amssymb}
\usepackage[mathscr]{eucal}
\usepackage{array,pinlabel,enumerate,color}
\usepackage[small,bf,margin=5pt]{caption}

\theoremstyle{plain}
\newtheorem{theorem}{Theorem}[section]
\newtheorem*{theorem*}{Theorem}

\newtheorem{corollary}[theorem]{Corollary}
\newtheorem*{corollary*}{Corollary}
\newtheorem{proposition}[theorem]{Proposition}
\newtheorem*{conjecture*}{Conjecture}
\theoremstyle{definition}\newtheorem{Def}[theorem]{Definition}
\theoremstyle{remark}
\newtheorem{example}[theorem]{Example}

\newcommand*{\N}{\ensuremath{\mathbf N}}
\newcommand*{\R}{\ensuremath{\mathbf R}}
\newcommand*{\Z}{\ensuremath{\mathbf Z}}
\newcommand*{\spin}{\ensuremath{\mathrm{Spin}^c}}

\DeclareMathOperator{\bip}{Bip}
\DeclareMathOperator{\conv}{Conv}
\DeclareMathOperator{\inter}{int}

\definecolor{gray}{gray}{.8}

\begin{document}

\title{Sutured Floer homology and hypergraphs}

\date{}

\author{Andr\'as Juh\'asz}
\thanks{AJ was supported by a Royal Society Research Fellowship and OTKA grant NK81203}
\address{Cambridge University}
\email{A.Juhasz@dpmms.cam.ac.uk}

\author{Tam\'as K\'alm\'an}
\thanks{TK was supported by the Japan Society for the Promotion of Science (JSPS) Grant-in-Aid for Young Scientists B, no.\ 21740041.}
\address{Tokyo Institute of Technology}
\email{kalman@math.titech.ac.jp}
\urladdr{www.math.titech.ac.jp/\char126 kalman}

\author{Jacob Rasmussen}
\address{Cambridge University}
\email{J.Rasmussen@dpmms.cam.ac.uk}


\begin{abstract}
By applying Seifert's algorithm to a special alternating diagram of a link $L,$ one obtains a Seifert surface $F$ of $L$.
We show that the support of the sutured Floer homology of the sutured manifold complementary to $F$ is affine isomorphic to
the set of lattice points given as hypertrees in a certain hypergraph that is naturally associated to the diagram.
This implies
that the Floer groups in question are supported in a set of $\spin$ structures that are the integer lattice points of a convex polytope.
This property has an immediate extension to Seifert surfaces arising from homogeneous link diagrams (including all alternating and
positive diagrams). 

In another direction, together with work in progress of the second author and others, our correspondence suggests a method for computing the ``top'' coefficients
of the HOMFLY polynomial of a special alternating link from the sutured Floer homology of a Seifert surface complement for a certain dual link.
\end{abstract}

\maketitle

\section{Introduction}

In this paper, we report on an unexpected coincidence between two
sets of integer lattice points that appear in the study of
alternating links and their Seifert surfaces.  The first set derives from the first
author's sutured Floer homology theory \cite{sfh}. Given a link $L
\subset S^3$ with Seifert surface $F$, we can associate to it the
sutured Floer homology group $SFH(S^3-F,L),$ where $(S^3-F,L)$ is the sutured manifold
complementary to $F$. This group decomposes
as a direct sum over relative $\spin$ structures on $(S^3-F,L),$ namely,
\[SFH(S^3-F,L) = \bigoplus_{\mathfrak{s} \in \spin(S^3-F,L)} SFH(S^3-F,L,\mathfrak{s}).\]
The set \(\spin(S^3-F,L)\) of relative $\spin$ structures can be thought of
as an affine copy of \(H_1(S^3-F;\Z)\). The \emph{support}
\[S(F)=\{\,\mathfrak{s} \in \spin(S^3-F,L) \mid SFH(S^3-F,L,\mathfrak{s}) \neq 0\,\}\]
carries a lot of interesting geometric information about \(F\), see \cite{sfhdecomp, sfhseifert, sfhpolytope}.

{\it A priori}, \(SFH(S^3-F,L,\mathfrak{s}) \) is a \(\Z/2\Z\)
graded group. However, if \(L\) is a non-split alternating link and \(F\) is an arbitrary
minimal genus Seifert surface of $L,$ then by \cite[Corollary 6.11]{decat}, each
group \(SFH(S^3-F,L,\mathfrak{s})\) is either trivial or isomorphic to
\(\Z,\) and all nontrivial groups have the same \(\Z/2\Z\)
grading. 
Thus, in this case, \(SFH(S^3-F,L)\) is determined by \(S(F)\). Furthermore, we see that $S(F)$ may be viewed as the Euler characteristic (function) $\chi(\mathfrak s)=\chi(SFH(S^3-F,L,\mathfrak{s}))$. Indeed, up to an overall sign, $\chi$ agrees with the indicator function of $S(F)$.


The second set was introduced by A.\ Postnikov \cite{post} and (independently, but later) by the second author \cite{hipertutte}, who based on it  a new theory of polynomial invariants of bipartite graphs and
hypergraphs. Its elements are generalized spanning trees called
hypertrees. More precisely, if $\mathscr H=(V,E)$ is a hypergraph with vertex set $V$ and hyperedge set $E$ (i.e., $V$ is a finite set and $E$ is a finite multiset of subsets of $V$), then an obvious two-to-one correspondence associates to it the bipartite graph $\bip\mathscr H$ with color classes $E$ and $V$. Now a hypertree in $\mathscr H$ is a vector $\mathbf f\colon E\to\N=\{\,0,1,2,\ldots\,\}$ so that $\bip\mathscr H$ has a spanning tree with valence $\mathbf f(e)+1$ at every $e\in E$. The set of all hypertrees in $\mathscr H$ will be denoted by\footnote{This slightly deviates from notation in \cite{hipertutte}, where $Q_\mathscr H$ stood for the convex hull of the set of hypertrees.} $Q_\mathscr H=Q_{(V,E)}$.

To place the support $S$ of $SFH$ and the set of hypertrees $Q$ in the same context, recall that to
any plane graph \(G\) we may associate an alternating link \(L_G\) and a surface \(F_G\). The latter is
obtained by taking a regular neighborhood of \(G\)
and inserting a
twist over each edge; \(L_G\) is its boundary. This is also known as the median construction on $G$. When $G$ is bipartite, $L_G$ is naturally oriented so that $F_G$ is its Seifert surface. The class of oriented links that arise as $L_G$ for some bipartite plane graph $G$ is called \emph{special alternating}, and it is known that $F_G$ is of minimal genus among Seifert surfaces of $L_G$. (Here, we allow multiple edges in $G$.) Let $R$ denote the set of connected components (regions) of $S^2-G$ and let $(E,R)$ be the hypergraph where $r \in R$ contains all elements of $E$ that lie along~$\partial r$. Define $(V,R)$ similarly. (We will soon shift the meaning of $r$ from region to a point marking the region, and thus $R$ will mean the set of those points.)
Note that $S^3-F_G$ is a handlebody of genus $|R|-1$.
Our main result is the following.

\begin{theorem}\label{thm:main}
Let \(G\) be a plane bipartite graph with color classes $E,V$ and regions $R$.  Then  $$S(F_G)\cong Q_{(E,R)}\cong-Q_{(V,R)}.$$
\end{theorem}

Here, the first isomorphism means that the $(|R|-1)$-dimensional affine lattice $\spin(S^3-F_G,L) \cong H_1(S^3-F_G,\Z)$ containing $S(F_G)$ has a $\Z$-affine identification with a certain hyperplane in $\Z^R$ so that the image of $S(F_G)$ is $Q_{(E,R)}$. The second isomorphism means that the two sets are translates of each other in $\Z^R$; this is quoted from \cite[Theorem 8.3]{hipertutte} to show that it hardly matters whether we pick $E$ or $V$ to play the role of the vertices in our hypergraph.

Both sets that appear in Theorem \ref{thm:main} can be viewed as
multi-variable polynomials and have representations in determinant
form. The theorem is proved by showing that a suitable
sequence of elimination steps transforms the `enhanced adjacency
matrix' of \cite{hipertutte} (which is a large matrix with simple monomial
entries) into the Turaev torsion (a smaller matrix with entries
obtained as Fox derivatives) that appears in \cite{decat}. It would be
very interesting to see if the relationship between the two matrices
is a special case of some more general phenomenon.

The statement of Theorem \ref{thm:main} may seem obscure, but in fact it has some far-reaching consequences.
Firstly, as $Q_{(E,R)}$ is easy to compute, and since the support $S(F_G)$ completely determines $SFH(S^3-F_G,L_G),$
we obtain a simple algorithm for determining the $SFH$ of complements of Seifert surfaces
given by the median construction on a bipartite graph. (Proposition \ref{pro:lengthen} illustrates the usefulness of this approach.) According to Hirasawa--Sakuma \cite{hisa} and Banks \cite{banks}, every minimal genus Seifert surface of a non-split, prime, special alternating link arises in this manner.

Together with the $\spin$-refined Murasugi and connected sum formulas obtained by using \cite[Proposition 5.4]{sfhpolytope} in the proofs of \cite[Corollary 8.8]{sfhdecomp} and \cite[Proposition 9.15]{sfh}, we can
also easily compute $SFH$ for Seifert surfaces obtained by applying Seifert's algorithm to homogeneous link diagrams \cite{cromwell} (which
include all alternating and positive diagrams). We will call such surfaces \emph{standard}.
Indeed, all standard surfaces are Murasugi sums and distant unions of standard Seifert surfaces of non-split, prime, special alternating links. When taking the Murasugi sum or distant union $F$ of the Seifert surfaces
$F_1$ and $F_2$ of the links $L_1$ and $L_2,$ respectively, and $L = \partial F,$ then
$$\spin(S^3-F,L) \cong \spin(S^3-F_1,L_1) \times \spin(S^3-F_2,L_2).$$
In case of a Murasugi sum,
$$SFH(S^3-F,L,(\mathfrak{s}_1,\mathfrak{s}_2)) \cong SFH(S^3-F_1,L_1,\mathfrak{s}_1) \otimes SFH(S^3 - F_2,L_2,\mathfrak{s}_2)$$
for every $\mathfrak{s}_i \in \spin(S^3-F_i,L_i)$, $i = 1,2$. In particular, notice that
\begin{equation} \label{eqn:Murasugi}
S(F) \cong S(F_1) \times S(F_2).
\end{equation}
On the other hand, if $F$ is the distant union of $F_1$ and $F_2,$ then using the connected sum formula for sutured manifolds,
$$SFH(S^3-F,L,(\mathfrak{s}_1,\mathfrak{s}_2)) \cong SFH(S^3-F_1,L_1,\mathfrak{s}_1) \otimes SFH(S^3 - F_2,L_2,\mathfrak{s}_2) \otimes \Z^2.$$
We again have (\ref{eqn:Murasugi}) for the supports, though $(S^3-F,L)$ is no longer a sutured $L$-space as the generators
of the $\Z^2$ factor lie in different $\Z/2\Z$-gradings.

As our second application, we derive
some previously unknown properties of the \(SFH\) of standard Seifert surfaces of homogeneous links.
We say that \(P \subset \Z^d\) is \emph{convex} if \(P \) is the intersection of \(\Z^d\) with the convex hull (in $\R^d$) of \(P\). One of the first non-trivial facts about $Q_\mathscr H$ is that it is convex, see \cite[Theorem 3.4]{hipertutte}. Combining this with (\ref{eqn:Murasugi}), we obtain the following.
\begin{corollary}
\label{Cor:convex}
If \(F\) is a standard Seifert surface of a homogeneous link L,
then \(S(F)\) is convex.
\end{corollary}
In contrast, if \(F\) is a Seifert surface obtained by applying
Seifert's algorithm to a non-homogeneous link \(L\), the support \(S(F)\) need not
be convex, even if the sutured Floer homology is supported in a single
\(\Z/2\Z\) grading. (For examples, see the computations for three-strand
pretzel knots in the last section of \cite{decat}.)


An unpublished theorem of
Richard Webb \cite{Webb} claims that all standard Seifert surfaces $F$ of
an alternating link have isomorphic $S(F)$. Thus if
\(L\) is a special alternating link, \(SFH(S^3-F,L)\), together with the $\spin$-grading,
does not depend on the choice of minimal genus Seifert surface $F$.
On the other hand, Altman \cite{irida} gives an example of
a knot \(K\) with distinct minimal genus Seifert
surfaces $F_1$ and $F_2$ for which $S(F_1)$ is convex and
$S(F_2)$ is not.
It would be interesting to know whether an
alternating link can have non-standard Seifert surfaces whose $SFH$
differs from that of the standard ones, or whether a homogeneous link can have a non-standard Seifert surface $F$ so that $S(F)$ is non-convex. 

Apart from convexity, there is another important context into which Theorem~\ref{thm:main} fits. Hypertrees were defined as part of a project aimed at finding
spanning tree models of some orientation-sensitive link invariants,
such as the HOMFLY polynomial.
Namely, by counting hypertrees appropriately, from $Q_\mathscr H$ we may read off a one-variable polynomial invariant $I_\mathscr H(\xi)$ of the hypergraph $\mathscr H$ \cite{hipertutte}. This so-called interior polynomial generalizes the partial evaluation $T(x,1)$ of the Tutte polynomial $T(x,y)$ of ordinary graphs. Then, extending a result of Jaeger \cite{jaeger}, the second author proposed the following.

\begin{conjecture*}
With notation as above, let $P_{L_G}(v,z)$ denote the HOMFLY polynomial of the link $L_G$. Then the part of $P_{L_G}$ which (after substituting $v=1$) becomes the leading term in the Alexander--Conway polynomial $\nabla_{L_G}(z)=P_{L_G}(1,z)$ is equal to $(vz)^{|R|-1}I_{(V,E)}(v^2)$.
\end{conjecture*}

As $E$ and $V$ play symmetric roles for $L_G$, this leads us to expect that $I_{(V,E)}=I_{(E,V)}$. Postnikov and the second author conjecture that this is indeed the case for any connected bipartite graph $G$; and furthermore, that the two polynomials also coincide with the so-called $h$-vector of some (any) triangulation of another polytope (called root polytope) derived from $G$. Then, it has already been proved~\cite{km} that this $h$-vector is equivalent to the relevant part of the HOMFLY polynomial. Therefore, the only missing step is purely discrete mathematical.

Now, if the Conjecture holds true, then
the part (hereafter called the \emph{top}) of $P_{L_G}(v,z)$
that realizes its maximum $z$-degree can be read off from either $Q_{(V,E)}$ or $Q_{(E,V)}$. According to Theorem \ref{thm:main}, these are equivalent to $S(F_{G_E})$ and $S(F_{G_V})$, respectively, where $G_E=\bip(V,R)=\bip(R,V)$ and $G_V=\bip(E,R)=\bip(R,E)$ are two plane bipartite graphs closely related to our original graph $G=G_R=\bip(V,E)=\bip(E,V)$. The three graphs together form a structure called a trinity. See Figure \ref{fig:links} for an example.

In other words, the Conjecture implies that the top of $P_{L_G}$ can be computed from Floer homology. The alert reader will notice that neither $F_{G_E}$ nor $F_{G_V}$ is a Seifert surface for $L_G=L_{G_R}$. It is also among our goals to derive the same HOMFLY coefficients from $S(F_{G_R})$. There are plenty of indications that this should be possible, including the following consequence of Theorem \ref{thm:main} and results in \cite{hipertutte}.

\begin{corollary}\label{cor:samesize}
With the plane bipartite graphs $G_R$, $G_E$, $G_V$ and associated surfaces defined as above, we have $|S(F_{G_R})|=|S(F_{G_E})|=|S(F_{G_V})|$ for the corresponding supports.
\end{corollary}

This is true despite the fact that the three sets are the lattice points in three polytopes of different dimensions. Corollary \ref{cor:samesize} is inspired by, and can be viewed as an extension of, Tutte's Tree Trinity Theorem \cite{tutte}. It is also closely related to Postnikov's duality result $|Q_{(V,E)}|=|Q_{(E,V)}|,$ see \cite{post}.

Murasugi and Przytycki \cite{mprz} showed that under `star product' of links, which is essentially Murasugi sum, the top of the HOMFLY polynomial behaves multiplicatively. Thus, for any homogeneous link, the Conjecture implies that coefficients along the top of the HOMFLY polynomial can be derived from Floer homology groups.

The paper is organized as follows. In sections \ref{sec:hyp} and \ref{sec:sfh} we review the necessary results about hypergraphs and Floer homology, respectively. Theorem \ref{thm:main}, along with some corollaries, is established in section \ref{sec:main}. Section \ref{sec:ex} contains a detailed sample calculation.


\pagebreak

\section{Hypergraph invariants}\label{sec:hyp}

In \cite{hipertutte}, the second author developed a new theory of hypergraph invariants by generalizing the spanning tree polytope and (partially) the Tutte polynomial to that context. Presently, we will only require the first item, the so called hypertree polytope. This was first defined by A.\ Postnikov \cite{post}, who called it the trimmed generalized permutohedron. (In this paper, instead of polytopes, we will only deal with their integer lattice points.) The relevant definitions and facts are as follows.

A \emph{hypergraph} is a pair $\mathscr H=(V,E)$, where $V$ is a finite set and $E$ is a finite multiset of non-empty subsets of $V$. Elements of $V$ are called \emph{vertices} and the elements of $E$ are the \emph{hyperedges}. Bipartite graphs and hypergraphs can be represented by the same type of picture. We formalize this idea in a two-to-one correspondence where each hypergraph $\mathscr H$ has an \emph{associated bipartite graph} $\bip\mathscr H$ in which $V$ and $E$ become the two color classes and $e\in E$ and $v\in V$ are connected if and only if $v\in e$. We will always assume that $\bip\mathscr H$ is connected. Conversely, a connected bipartite graph $G$ \emph{induces} 
two hypergraphs which we call abstract duals. In other words, the abstract dual $\overline{\mathscr H}=(E,V)$ of $\mathscr H$ is defined by interchanging the roles of its vertices and hyperedges.

\begin{Def}
Let $\mathscr H=(V,E)$ be a hypergraph so that its associated bipartite graph $\bip\mathscr H$ is connected. By a \emph{hypertree} in $\mathscr H$ we mean a function (vector) $\mathbf f\colon E\to\N=\{\,0,1,\ldots\,\}$ so that a spanning tree of $\bip\mathscr H$ can be found which has valence $\mathbf f(e)+1$ at each $e\in E$.
The set of all hypertrees in $\mathscr H$ will be denoted with $Q_\mathscr H$.
\end{Def}

An elementary observation is that all hypertrees lie on the hyperplane
\[\Pi=\left\{\,\mathbf f\colon E\to\N\biggm|\sum_{e\in E}\mathbf f(e)=|V|-1\,\right\}\subset\R^E.\]
Furthermore, it turns out that hypertrees are exactly the integer lattice points of a
polytope
which is cut out of $\Pi$ by inequalities of the form $\sum_{e\in E'}\mathbf f(e)\le\mu(E')$, where $\varnothing\ne E'\subsetneqq E$. Here $\mu$ is defined as follows\footnote{If we let $\mu(\varnothing)=0$ and $\mu(E)=|V|-1$, we obtain a non-decreasing submodular set function.
}. Let $\bip\mathscr H\big|_{E'}$ be the bipartite graph with color classes $E'\subset E$ and $\bigcup E'\subset V$, and edges inherited from $\bip\mathscr H$. Let us denote the number of its connected components by $c(E')$. Then, set $\mu(E')=|\bigcup E'|-c(E')$.

\begin{figure}[htbp]
\labellist
\footnotesize
\pinlabel $r_0$ at 366 98
\pinlabel $r_1$ at 259 422
\pinlabel $r_2$ at 601 368
\pinlabel $r_3$ at 528 633
\pinlabel $v_0$ at 276 602
\pinlabel $v_1$ at 389 462
\pinlabel $v_2$ at 352 314
\pinlabel $v_3$ at 690 206
\pinlabel $v_4$ at 623 548
\pinlabel $e_0$ at 834 529
\pinlabel $e_1$ at 420 529
\pinlabel $e_2$ at 366 390
\pinlabel $e_3$ at 347 242
\pinlabel $t_0$ at 200 750
\pinlabel $t_1$ at 400 595
\pinlabel $t_2$ at 350 490
\pinlabel $t_3$ at 175 410
\pinlabel $t_4$ at 320 365
\pinlabel $t_5$ at 450 150
\pinlabel $t_6$ at 410 282
\pinlabel $t_7$ at 450 410
\pinlabel $t_8$ at 570 490
\pinlabel $t_9$ at 670 595
\pinlabel $t_{10}$ at 740 340
\endlabellist
   \centering
   \includegraphics[width=\linewidth]{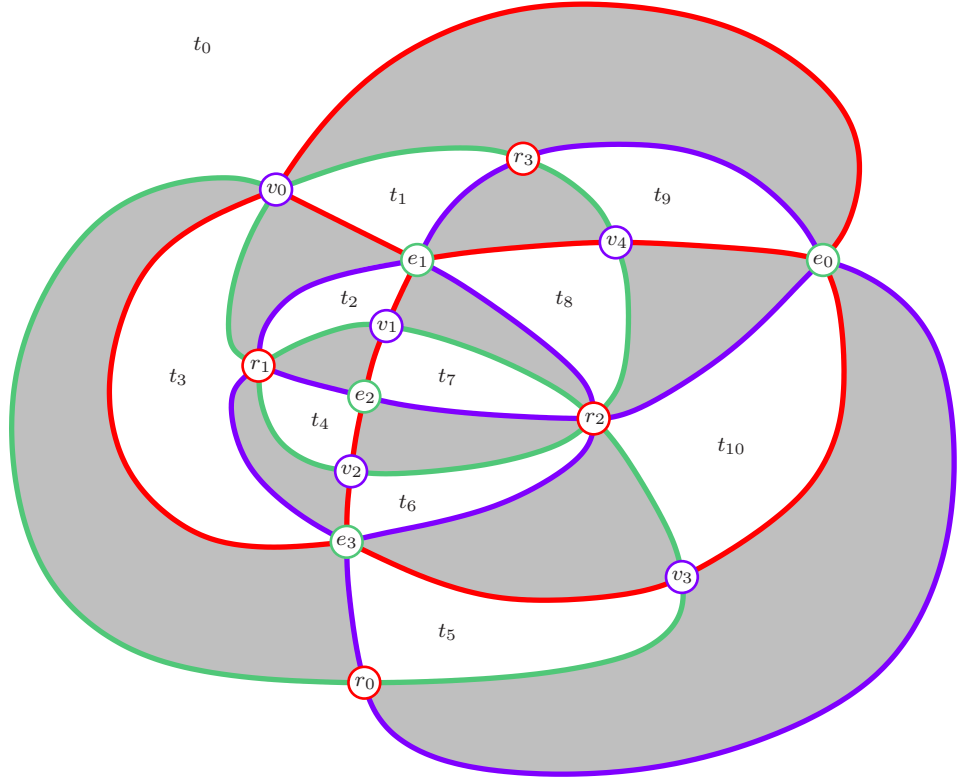}
   \caption{A trinity of plane bipartite graphs.}
   \label{fig:trinity}
\end{figure}

For knot theoretical considerations, a case of particular interest is when $\bip\mathscr H$ is a planar graph, i.e., it comes with an embedding in $S^2$. Such graphs always occur in sets of three, forming a structure called a trinity. Trinities, first investigated by Tutte \cite{tutte}, have many equivalent descriptions; we choose to define them as triangulations of $S^2$ with a proper three-coloring of the $0$-cells (which we will generally call `points,' with the understanding that some of them may become `vertices' in a graph or hypergraph). See Figure \ref{fig:trinity} for an example.

We will use the names red, emerald and violet for the colors,
and denote the respective sets of points by $R$, $E$ and $V$. Let us color each edge in the triangulation with the color that does not occur among its ends. Then $E$ and $V$ together with the red edges form a bipartite graph that we will call the \emph{red graph} and denote by~$G_R$. Each region (i.e., connected component of the complement) of the red graph contains a unique red point. Likewise, the \emph{emerald graph} $G_E$ has red and violet points, emerald edges, and regions marked with emerald points. Finally, the \emph{violet graph} contains $R$ and $E$ as vertices, violet edges, and a violet point in each of its regions.

Any one of the six hypergraphs induced by $G_R$, $G_E$, or $G_V$ generates the other five (and only those) using two operations. One is abstract duality and the other one is called planar duality. In the latter, hyperedges are retained, but the vertices of the hypergraph are replaced with (points placed in) the regions of its associated plane bipartite graph; the containment relation is defined by adjacency.

In \cite[Theorem 8.3]{hipertutte}, it is established that the hypertree polytopes of planar dual hypergraphs are centrally symmetric. Thus in any trinity, we have
$Q_{(V,E)}\cong-Q_{(R,E)}\subset\R^E$, $Q_{(R,V)}\cong-Q_{(E,V)}\subset\R^V$, and
\begin{equation}\label{symm}
Q_{(E,R)}\cong-Q_{(V,R)}\subset\R^R.
\end{equation}

Notice that each triangle of a trinity is adjacent to exactly one edge and one point of each color. Compared to the orientation of the sphere, the cyclic order of the colors around each triangle may be positive or negative. If two triangles share an edge, these orientations are opposite.
Hence the triangles have a black and white checkerboard coloring according to orientation, cf.\ Figure \ref{fig:trinity}. This coloring can be used to orient the dual graphs $G_R^*$, $G_E^*$, and $G_V^*$ (dual in the classical, not in the hypergraph sense). Namely, each edge in these graphs cuts through a black and a white triangle, and we orient it so that its tail end is in black territory.

Next, note that the sets of red edges, emerald edges, violet edges, white triangles, and black triangles all have the same cardinality $n$. In particular, adjacency defines natural bijections between white triangles and edges of each color. Now, if we apply Euler's formula to the trinity, we get $|R|+|E|+|V|-3n+2n=2$, i.e,
\begin{equation}\label{eq:euler}
\text{ the total number of points exceeds that of the white triangles by }2.
\end{equation}

Let us distinguish one white triangle as the \emph{outer} one and label it with $t_0$. We will label the red, emerald, and violet points adjacent to $t_0$ with $r_0$, $e_0$, and $v_0$, respectively, and call them the \emph{roots}. According to \eqref{eq:euler}, the non-outer white triangles and non-root points are equal in number, so they form a square adjacency matrix $M$. We index the rows of $M$ by points and its columns by triangles. This matrix was first investigated by Berman \cite{berman} who showed that when we compute the determinant $\det M$, all non-zero expansion terms appear with the same sign.
We call these terms the \emph{Tutte matchings} of the trinity.

A Tutte matching is indeed a complete matching of non-root points to adjacent non-outer white triangles. Equivalently, it can be thought of as a collection of edges from $G_R^*$, $G_E^*$ and $G_V^*$. Namely, whenever a triangle and a point are matched, we choose the dual edge that cuts through the triangle and ends at the point. It turns out that such a collection is the union of three so called spanning arborescences of $G_R^*$, $G_E^*$ and $G_V^*$, respectively, rooted at the respective points adjacent to $t_0$. (Given an arbitrary vertex $q$, called root, in a directed graph $D$, a \emph{spanning arborescence} is a spanning tree of $D$ which contains a directed path from $q$ to any other vertex $p\ne q$ of $D$.) Each spanning arborescence occurs in relation to exactly one Tutte matching. This fact is the basis of the proof of Tutte's Tree Trinity Theorem \cite{tutte}, i.e., of the claim that $G_R^*$, $G_E^*$ and $G_V^*$ have the same number of spanning arborescences. (These counts do not depend on the position of the roots either, and of course they are equal to $|\det M|$.)

\begin{Def}
If a non-outer white triangle $t_i$ is adjacent to the red point $r_j$ and the non-root emerald point $e_k$, then in $M$, at the intersection of row $e_k$ and column $t_i$, let us change the entry $1$ to $r_j$. After it becomes a matrix entry, we will think of $r_j$ as an indeterminate associated with the original point. Call the resulting matrix the \emph{enhanced adjacency matrix} and denote it with $M_{r\to e}$.
\end{Def}

By varying the colors, altogether six such matrices can be associated to a trinity. According to \cite[Theorem 10.5]{hipertutte}, we may use enhanced adjacency matrices to compute the set of hypertrees in any hypergraph $\mathscr H$ as long as $\bip\mathscr H$ is planar.

\begin{proposition}\label{pro:det}
In any trinity, for the set of hypertrees in the hypergraph with emerald vertices and red hyperedges, we have
\[Q_{(E,R)}
=\det M_{r\to e}\]
in the following sense. The determinant on the right hand side is a sum of monomials in the indeterminates $r\in R$. Either each monomial has coefficient $+1,$ or each has coefficient $-1$. If we write the exponents in the monomials as vectors, the set we obtain is exactly the left hand side.
\end{proposition}

\begin{example}\label{ex:det}
The adjacency matrix associated with the trinity of Figure \ref{fig:trinity} is
\[M=\bordermatrix
{~&\text{\footnotesize $t_1$}&\text{\footnotesize $t_2$}&\text{\footnotesize $t_3$}&\text{\footnotesize $t_4$}&\text{\footnotesize $t_5$}&\text{\footnotesize $t_6$}&\text{\footnotesize $t_7$}&\text{\footnotesize $t_8$}&\text{\footnotesize $t_9$}&\text{\footnotesize $t_{10}$}\cr
\text{\footnotesize $r_1$}&0&1&1&1&0&0&0&0&0&0\cr
\text{\footnotesize $r_2$}&0&0&0&0&0&1&1&1&0&1\cr
\text{\footnotesize $r_3$}&1&0&0&0&0&0&0&0&1&0\cr
\text{\footnotesize $e_1$}&1&1&0&0&0&0&0&1&0&0\cr
\text{\footnotesize $e_2$}&0&0&0&1&0&0&1&0&0&0\cr
\text{\footnotesize $e_3$}&0&0&1&0&1&1&0&0&0&0\cr
\text{\footnotesize $v_1$}&0&1&0&0&0&0&1&0&0&0\cr
\text{\footnotesize $v_2$}&0&0&0&1&0&1&0&0&0&0\cr
\text{\footnotesize $v_3$}&0&0&0&0&1&0&0&0&0&1\cr
\text{\footnotesize $v_4$}&0&0&0&0&0&0&0&1&1&0\cr}.\]
Its determinant is $11$.
One of its enhancements is as follows.
\begin{equation}\label{adjmatrix}
M_{r\to e}=\bordermatrix
{~&\text{\footnotesize $t_1$}&\text{\footnotesize $t_2$}&\text{\footnotesize $t_3$}&\text{\footnotesize $t_4$}&\text{\footnotesize $t_5$}&\text{\footnotesize $t_6$}&\text{\footnotesize $t_7$}&\text{\footnotesize $t_8$}&\text{\footnotesize $t_9$}&\text{\footnotesize $t_{10}$}\cr
\text{\footnotesize $r_1$}&0&1&1&1&0&0&0&0&0&0\cr
\text{\footnotesize $r_2$}&0&0&0&0&0&1&1&1&0&1\cr
\text{\footnotesize $r_3$}&1&0&0&0&0&0&0&0&1&0\cr
\text{\footnotesize $e_1$}&r_3&r_1&0&0&0&0&0&r_2&0&0\cr
\text{\footnotesize $e_2$}&0&0&0&r_1&0&0&r_2&0&0&0\cr
\text{\footnotesize $e_3$}&0&0&r_1&0&r_0&r_2&0&0&0&0\cr
\text{\footnotesize $v_1$}&0&1&0&0&0&0&1&0&0&0\cr
\text{\footnotesize $v_2$}&0&0&0&1&0&1&0&0&0&0\cr
\text{\footnotesize $v_3$}&0&0&0&0&1&0&0&0&0&1\cr
\text{\footnotesize $v_4$}&0&0&0&0&0&0&0&1&1&0\cr}.
\end{equation}
The determinant of the latter is $r_0 r_1^2 + r_0 r_1 r_2 + r_1^2 r_2 + r_0 r_2^2 + r_1 r_2^2 + r_2^3 + r_0 r_1 r_3 +
 r_1^2 r_3 + r_0 r_2 r_3 + r_1 r_2 r_3 + r_2^2 r_3$. If we interpret this formula as in Proposition \ref{pro:det}, we obtain the picture in Figure \ref{fig:polytope} for the set $Q_{(E,R)}$ of hypertrees.
\begin{figure}[htbp]
\unitlength 7pt
\begin{picture}(24,27)
\color{gray}
\put(18,3){\line(0,1){21}}
\put(18,3){\line(2,3){2}}
\put(18,3){\line(-2,1){18}}
\put(18,24){\line(1,-2){4}}
\put(18,24){\line(-3,-2){18}}
\put(0,12){\line(1,0){11}}
\color{black}
\thicklines
\put(18,10){\line(-2,1){6}}
\put(16,12){\line(1,0){1}}
\put(24,12){\line(-1,0){4}}
\put(24,12){\line(-1,2){2}}
\put(18,10){\line(2,3){4}}
\put(12,6){\line(1,0){5}}
\put(12,6){\line(0,1){7}}
\put(12,6){\line(3,2){6}}
\put(12,6){\line(2,3){3}}
\put(20,6){\line(-1,0){1}}
\put(20,6){\line(2,3){4}}
\put(20,6){\line(-1,2){2}}
\put(14,16){\line(1,0){3}}
\put(14,16){\line(-2,-3){2}}
\put(14,16){\line(1,-2){2}}
\put(22,16){\line(-1,0){3}}
\thinlines
\put(24,12){\circle*{1}}
\put(18,10){\circle*{1}}
\put(16,12){\circle*{1}}
\put(12,13){\circle*{1}}
\put(20,13){\circle*{1}}
\put(12,6){\circle*{1}}
\put(20,6){\circle*{1}}
\put(22,9){\circle*{1}}
\put(22,16){\circle*{1}}
\put(14,9){\circle*{1}}
\put(14,16){\circle*{1}}
\put(-2,12){$r_0^3$}
\put(18,1){$r_1^3$}
\put(25,12){$r_2^3$}
\put(18,25){$r_3^3$}
\end{picture}
\caption{The eleven hypertrees in a hypergraph and their convex hull.}
\label{fig:polytope}
\end{figure}
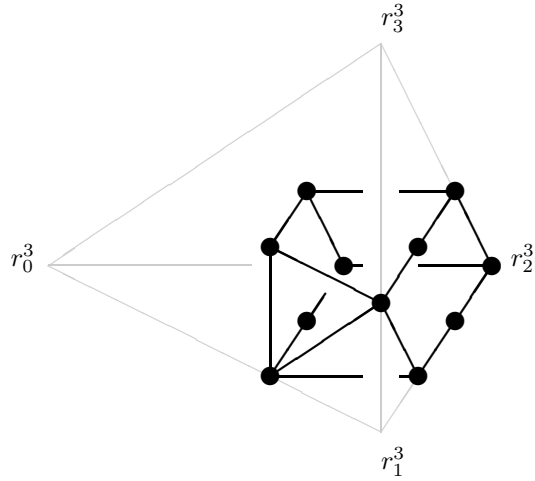
\end{example}


\section{Decategorified sutured Floer homology}\label{sec:sfh}

Sutured Floer homology \cite{sfh} is an invariant of balanced sutured
manifolds. A \emph{sutured manifold}  $(M,\gamma)$ is an oriented three-manifold
\(M\), together with an oriented, null-homologous  multicurve \(\gamma\)
(the \emph{sutures}) on \(\partial M,\) see \cite{gabai}. Since \(\gamma\)
is null-homologous, it divides \(\partial M - \gamma\) into two
pieces \(R_+\) and \(R_-\), with the property that \(R_+\) lies on one
side of each curve in \(\gamma\), and \(R_-\) lies on the other.
(Note
that \(R_\pm\) need not be connected.)
The condition that \((M,
\gamma)\) is \emph{balanced} means that
\begin{enumerate}[(i)]
\item \(\chi(R_+)=\chi(R_-),\)
\item every component of \(\partial M\) contains at least one suture and
\item $M$ has no closed components.
\end{enumerate}

The sutured manifolds we consider in this paper all arise from the
following construction.
Suppose that \(L \subset S^3\) is an oriented link, and that \(F\) is
a Seifert surface for \(L\). Let us thicken \(F\) slightly to \(F
\times [\,-\varepsilon, \varepsilon\,] \subset S^3\), and let \(M_F\) be the
closure of its complement. The link \(L\) is a subset of \(\partial
M_F\), and the pair \((M_F,L)\) is a balanced sutured
manifold. Note that both subsurfaces \(R_\pm\) are homeomorphic to \(\inter F\). Thus in this context, we will prefer to use the notation $R_+=F_+=F\times\{\varepsilon\}$ and $R_-=F_-=F\times\{-\varepsilon\}$.

When \(L\) is an alternating link and \(F\) is a minimal genus Seifert
surface of $L,$ the sutured Floer homology groups \(SFH(M_F,L)\) are
completely determined by their Euler characteristic according to \cite[Corollary 6.6]{decat}. The main result of \cite{decat} identifies this Euler
characteristic with a certain Turaev torsion \(\tau(M_F,L)\). For our
purposes, \(\tau(M_F,L)\) is best viewed as an element of the group
ring \(\Z[H_1(M_F;\Z)]\) which is well defined up to multiplication by a
unit in the group ring.

Let now $G=G_R$ be the red graph of a trinity. As we saw in the previous section, any connected plane bipartite graph uniquely arises in this way. Namely, color the edges of $G$ in red and its two color classes in emerald and violet; place a red point in each region of $G$ and make the appropriate emerald and violet connections between these new points and the old ones. In what follows, let us assume that $G$ is $2$-connected, i.e., that it has no cut vertices.

\begin{figure}[htbp] 
   \centering
   \includegraphics[width=1in]{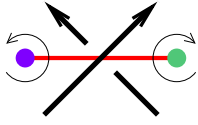}
   \caption{The two strands of $L_G$ crossing an edge of $G$.}
   \label{fig:crossing}
\end{figure}

We use the so called median construction to associate the alternating link $L_G$ to~$G$. I.e., along the boundary of each region of $G$, we first connect the midpoints of consecutive edges by disjoint simple curves. There are exactly two ways to specify over- and undercrossing information at the midpoints themselves so that the union of these curves becomes an alternating link. We use the convention of Figure \ref{fig:crossing}. We also place a positive spin (as in a small counterclockwise spinning top) at each violet point and a negative spin at each emerald point. The link $L_G$ inherits an orientation from the spins as shown in Figure \ref{fig:crossing} so that its diagram becomes positive. Seifert circles of the diagram correspond to the emerald and violet points so that they are not nested in one another. Hence $L_G$ is a special alternating link, and in fact, any positive special alternating link arises as the result of this construction. (Changing one or both of our choices to the opposite, i.e., mirroring and/or reversing $L_G$, does not make any major difference.)

\begin{figure}[htbp]
   \centering
   \includegraphics[width=\linewidth]{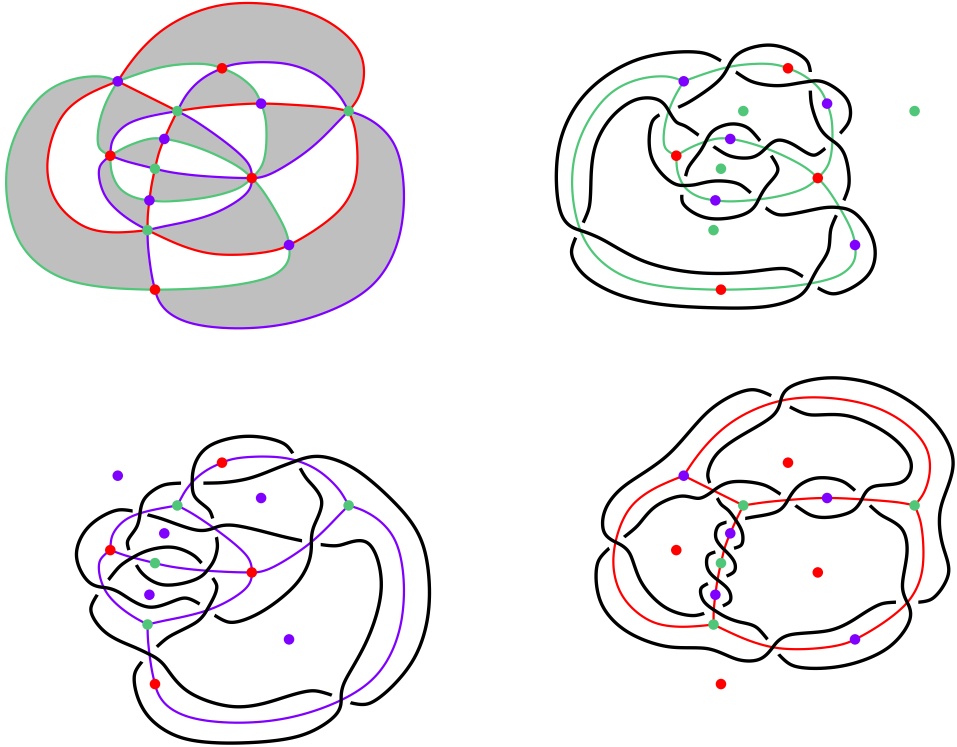}
   \caption{The bipartite graphs of a trinity and the associated special alternating links.}
   \label{fig:links}
\end{figure}

By repeating our construction for the other two colors, we end up associating three distinct special alternating links to our trinity. See Figure \ref{fig:links} for an example. These links are of course closely related; see, for instance, Corollary \ref{cor:alex}.

Our construction is such that $G$ is recovered as the
Seifert graph of $L_G$. Hence, if we apply Seifert's construction to
the diagram of $L_G,$ we obtain a Seifert surface $F_G$ which
deformation retracts onto $G$. Let  \(M_G=M_{F_G}\) be the complement
of this surface. It is easy to see that \(M_G\) is a handlebody of genus $|R|-1$.
We are interested in understanding the group \(SFH(M_G,L_G)\). Since
\(L_G\) is alternating, it suffices to study the torsion
\(\tau(M_G,L_G)\). We now recall from \cite{decat} how
 to compute \(\tau(M,\gamma)\)
in the case when  \(R_-\) is connected and \(M\) is a handlebody of
genus \(g\).

First, we choose a system of co-oriented compressing disks
\(A_1,\ldots, A_g \subset M\) so that when we cut \(M\) along the
\(A_i\), the result is a ball.
We also fix a basepoint \(p \in R_- \) disjoint from the \(A_i\). These data determine a specific isomorphism between
\(\pi_1(M,p)\) and the free group generated by elements
\(a_1,\ldots,a_g\). Namely, if \(c\) is a loop based at \(p\), we can
isotope \(c\) so it is transverse to the \(A_i\). Then the word
associated to \(c\) is obtained by traversing \(c\) and recording
either an \(a_i\) or an \(a_i^{-1}\) (depending on the sign of
intersection) each time we pass through \(A_i\).

Since \((M,\gamma)\) is balanced, \(\chi(R_-) =
\frac{1}{2}\chi(\partial M) = 1-g\). It follows that \(\pi_1(R_-,p)\) is a free
group on \(g\) generators. Choose a set of generators for \(\pi_1(R_-,p)\)  and consider
their images  \(W_1,\ldots,W_g \in \pi_1(M,p) \cong \langle
a_1,\ldots, a_g \rangle \).

\begin{proposition}\cite[Prop. 5.1]{decat} \label{pro:turaev}
If $M$ is a handlebody and $R_-$ is connected, then, up to multiplication by a unit, we have \[ \tau (M,\gamma) \sim \det (\left \vert
  {d_{a_i} W_j} \right \vert), \] where \(d_{a_i}\) denotes the Fox free
derivative with respect to \(a_i\) and \(|\cdot | \colon \Z[\pi_1(M)] \to
\Z[H_1(M)]\) is the abelianization.
\end{proposition}

\begin{example}\label{ex:paths}
Let us choose $G$ to be the red bipartite graph $G_R$ of Figure \ref{fig:links}. One way to compute the corresponding Turaev torsion is as follows. Please refer to Figure~\ref{fig:compute}. Our basepoint $p$ will be the push-off of $v_0$ on the side of $F$ facing the viewer. We fix loops on that side $F_-$ of the surface that go counterclockwise once around the boundary of the region marked with $r_1$, clockwise around $r_0$, and counterclockwise around $r_3$. (The second loop appears in the figure as a counterclockwise path around the outside contour of $G_R$.)

Our compressing disks correspond to the regions $r_1$, $r_2$, and $r_3$. We co-orient each toward the center of the sphere containing the diagram. (Our figures are drawn from a viewpoint which is outside of the sphere.) For the corresponding free group generators we will use the same symbols $r_j$.

With these conventions, the three loops that we described above yield the words
\[r_1 r_2^{-1} r_1 r_2^{-1} r_1 r_3^{-1},\quad r_1 r_0^{-1} r_2 r_0^{-1},\quad\text{and}\quad r_3 r_2^{-1} r_3 r_0^{-1},\]
respectively. Note that we also counted intersections with the region $r_0$, which did not contribute a compressing disk. In other words, $r_0=1$ but we will keep the symbol $r_0$ around as its presence helps the formalism in the proof of Theorem \ref{thm:main}. Finally, the Turaev torsion is
\begin{multline}\label{turaev}
\left|
\begin{array}{ccc}
1+r_1r_2^{-1}+r_1^2r_2^{-2}&1&0\\
-r_1r_2^{-1}-r_1^2r_2^{-2}&r_1r_0^{-1}&-r_3r_2^{-1}\\
-r_1^3r_2^{-2}r_3^{-1}&0&1+r_3r_2^{-1}
\end{array}
\right|=r_0^{-1}r_1+r_0^{-1}r_1^2r_2^{-1}+r_0^{-1}r_1^3r_2^{-2}\\
\\
+r_0^{-1}r_1r_2^{-1}r_3+r_0^{-1}r_1^2r_2^{-2}r_3+r_0^{-1}r_1^3r_2^{-3}r_3+r_1^3r_2^{-3}+r_1r_2^{-1}+r_1^2r_2^{-2}+r_1r_2^{-2}r_3+r_1^2r_2^{-3}r_3.
\end{multline}
The reader may check that \eqref{turaev} agrees with a monomial ($r_0^{-1}r_1r_2^{-3}$) times the second determinant of Example \ref{ex:det}. Our main result says that this is not a coincidence.
\end{example}

\section{The main result}\label{sec:main}

Due to Propositions \ref{pro:det} and \ref{pro:turaev}, we may re-state Theorem \ref{thm:main} in the following equivalent form.

\begin{theorem}
Let $G$ be a $2$-connected plane bipartite graph (with color classes $E$ and $V$) and let us realize it as the red graph $G=G_R$ of a trinity. Let the positive special alternating link $L_G$, Seifert surface $F_G$ and handlebody $M_G$ be associated to $G$ as above. Then, the Turaev torsion $\tau(M_G,L_G)$ is, up to multiplication by a monomial, equal to the determinant of the enhanced adjacency matrix $M_{r\to e}$.
\end{theorem}


By writing `equal' above, we assume various identifications through which both determinants become polynomials in indeterminates corresponding to the red points of the trinity. For $\det M_{r\to e}$ this occurs naturally; for $\tau(M_G,L_G)$ it is made explicit below. The proof is not deep but complicated. In the next section, we carry out a concrete computation based on the same method. The reader may wish to study that first or to read it parallel to the proof.

\begin{proof}
We are going to describe a very specific way of computing the Turaev torsion, along with an equally specific way of manipulating the enhanced adjacency matrix so that eventually the two determinants are shown to coincide up to a monomial factor. Let us fix an outer white triangle $t_0$, with adjacent points (roots) $r_0$, $e_0$, and $v_0$ as before. Other points will be labeled with $r_1,\ldots,r_{|R|-1}$, $e_1,\ldots,e_{|E|-1}$, and $v_1,\ldots,v_{|V|-1}$. All additional choices will be dictated by an arbitrarily fixed Tutte matching. I.e., we break the symmetry by singling out one expansion term in $\det M_{r\to e}$.

\emph{The Turaev torsion:} We compute $\tau(M_G,L_G)$ using the following procedure. First, let us choose compressing disks for $M_G$ the same way as in Example \ref{ex:paths}: Let $S^2$ denote the sphere containing $G,$ and let us fix an `outside viewpoint' (which is used to draw our diagrams) and a `center'. Connect these two points by a path through $r_0$ and thicken it slightly into a $3$-ball. We may think of $M_G$ as this $3$-ball with $|R|-1$ one-handles attached so that each handle passes through a non-root red point. For each $j=1,2,\ldots,|R|-1$, there is an obvious disk contained in $S^2$, centered at $r_j$, which is the co-core to one of these handles. We co-orient each disk toward the center of $S^2$. When the $j$'th compressing disk is identified with a free group generator, we will write $r_j$ for the latter as well.


Our loops in $F_-$ will be based at $(v_0,-\varepsilon)$. Here, we assume that the thickening of $F=F_G$ has been parametrized so that near the midpoint of the edge $\kappa$ of $t_0$ that connects $e_0$ and $v_0$, it is $F_-$ that faces the region 
$r_0$.

Next, fix a spanning arborescence $A$ in $G^*$ rooted at $r_0$. As we explained in Section \ref{sec:hyp}, $A$ assigns an adjacent non-outer white triangle to every non-root red point. We choose the labels of these triangles so that $t_j$ is assigned to
$r_j$ for all $j=1,2,\ldots,|R|-1$, i.e., that the unique edge of $A$ pointing toward $r_j$ reaches $r_j$ through $t_j$.

Let now $\Gamma\subset G$ denote the dual tree of $A$. Note that $\Gamma$ contains the edge $\kappa$ because $t_0$ is not matched to any point.
For all $j$, if we add the dual of the edge of $A$ pointing toward $r_j$ (suppose it has endpoints $e^*$ and $v^*$) to $\Gamma$, a unique cycle is created. We describe it as a sequence $\lambda_j$ of edges of $G$ that starts at $v_0$, follows the unique path in $\Gamma$ to $v^*$, contains the edge $v^* e^*$ as the next entry, finally returns from $e^*$ to $v_0$ along the unique path that exists in $\Gamma$. In Figure \ref{fig:snake} we indicate the two sequences of edges that, along with $v^*e^*$, make up the loop. Note that, in addition to the point $v_0$ appearing twice, it is possible for the first (leftmost) several edges to agree, in reverse order, with the last (rightmost) several edges. Other than these coincidences, the edges in $\lambda_j$, and hence the white triangles adjacent to them, are all different (in fact, it will sometimes be more useful to think of $\lambda_j$ as a sequence of triangles rather than of edges) but there may be repetitions among the red points adjacent to those white triangles. In particular, both $r_j$ and $r_0$ may occur several times.

\begin{figure}[htbp] 
\labellist
\tiny
\pinlabel $v_0$ at 25 143
\pinlabel $p_{l'}$ at 225 143
\pinlabel $t_l$ at 270 100
\pinlabel $x$ at 269 18
\pinlabel $p_l$ at 315 143
\pinlabel $p_{k'}$ at 423 143
\pinlabel $t_k$ at 468 200
\pinlabel $y$ at 467 269
\pinlabel $p_k$ at 513 143
\pinlabel $b$ at 665 271
\pinlabel $a$ at 756 18
\pinlabel $v^*$ at 801 145
\pinlabel $r_j$ at 846 268
\pinlabel $e^*$ at 891 145
\pinlabel $c$ at 935 18
\pinlabel $w$ at 980 144
\pinlabel $v_0$ at 1179 143
\endlabellist
   \centering
   \includegraphics[width=\linewidth]{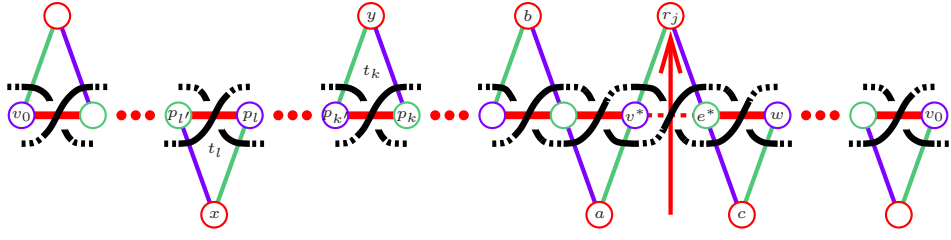}
   \caption{Schematic image of the loop $\lambda_j$ corresponding to the non-root red point $r_j$.}
   \label{fig:snake}
\end{figure}

It is easy to see that the loops $\lambda_j$ serve as free generators for $\pi_1(G,v_0)$ and hence their push-offs to $F_-$ freely generate $\pi_1(F_-)$. Let us proceed to the word $W_j$, $j=1,2,\ldots,|R|-1$, in the letters $r_0,r_1,\ldots,r_{|R|-1}$, that is derived from $\lambda_j$. (Here the letter $r_0$ will be used to keep track of intersections with the region $r_0$. We do this just to preserve symmetry. In the free group, $r_0=1$.) As $L_G$ is alternating, the exponents $+1$ and $-1$ alternate in $W_j$. Our conventions regarding the thickened surface have been set so that $F_-$ faces the observer near violet points, while near emerald points we see $F_+$. Hence in $W_j$ the first letter has exponent $+1$, the last has $-1$, and the edge $v^*e^*$ contributes $r_j^{+1}$. Highlighting that contribution, let us write $W_j=U_jr_jZ_j$.

When we consider the Fox derivatives of $W_j$ with respect to $r_1,\ldots,r_{|R|-1}$, we see that every letter in $W_j$ that is different from $r_0^{\pm1}$ contributes exactly one monomial to exactly one Fox derivative. Namely, if the red point $x$ appears on the right side of the path (cf.\ Figure \ref{fig:snake}) then the contribution is
\[-W^xx^{-1}\text{ to }\partial W_j/\partial x,\]
and a red point $y$ on the left side contributes
\[W^y\text{ to }\partial W_j/\partial y.\]
Here, $W^z$ denotes the part of $W_j$ that precedes the letter $z$. (This notation is sloppy because $z$ may appear several times in $W_j$ and thus there can be several different $W^z$'s; yet we hope that no confusion will arise.) Our task in the next part of the proof is to arrive at the same monomials in an entirely different way.

\emph{Manipulating the enhanced adjacency matrix:} We will carry out a sequence of elimination steps on $M_{r\to e}$ inspired by the proof of the Tree Trinity Theorem. We choose a \emph{deconstruction order} $<$ of the points $e_1,\ldots,e_{|E|-1},v_1,\ldots,v_{|V|-1}$, by which we mean an order so that the smallest point is a vertex of valence $1$ in $\Gamma$, then the second smallest is a vertex of valence $1$ in the tree that remains if we erase the first vertex and its adjacent edge from $\Gamma$, and so on. After removing all non-root vertices in this manner, only the edge $\kappa$ remains of $\Gamma$. For ease of exposition, we will assign a second label to non-root emerald and violet points so that
\[\{\,e_1,\ldots,e_{|E|-1},v_1,\ldots,v_{|V|-1}\,\}=\{\,p_{|R|}>p_{|R|+1}>\cdots>p_{n-2}>p_{n-1}\,\}.\]
Here, $n-1=|R|+|E|+|V|-3$ is the number of non-root points in the trinity.

The deconstruction order can be used to match non-outer white triangles to the non-root emerald and violet vertices in a one-to-one fashion. The non-root red vertices have already been matched to triangles by the arborescence $A$ and we will not use those triangles again, so that our construction will result in a Tutte matching. The tree $\Gamma$ contains an edge if and only if the white triangle adjacent to the edge has not been matched to a red point. For the emerald or violet point $p_{n-1}$ having valence $1$ in $\Gamma$ means that all but one of the adjacent white triangles have been matched to red points, so that there is only one choice left. We match that triangle to $p_{n-1}$. A similar unique choice of a triangle exists for $p_{n-2}$ and so on. We complete the labeling of the white triangles so that the process matches the point $p_i$ with the triangle $t_i$ for each $i=|R|,\ldots,n-1$.

The same matching can be described without reference to the deconstruction order as follows: given a non-root point $p_i$ in $\Gamma$, the tree has a unique edge adjacent to it which is such that its other endpoint is closer to $\kappa$ than $p_i$ is. The
white triangle adjacent to this edge is $t_i$.

From our assumptions on the labeling it follows that if we list the rows and columns of $M_{r\to e}$ in the orders $r_1,\ldots,r_{|R|-1},p_{|R|},\ldots,p_{n-1}$ and $t_1,\ldots,t_{n-1}$, respectively, then no zeros appear along the main diagonal and the bottom right $(n-|R|)\times(n-|R|)$ block is upper triangular. Indeed, if $i\ge|R|$, then the red edge of the triangle $t_i$ is in $\Gamma$ and $p_i$ is the endpoint of this edge which is farther away in $\Gamma$ from $\kappa$. Therefore the other emerald or violet point of $t_i$ is larger than $p_i$ in the deconstruction order.

We will use the $n-|R|$ diagonal entries of this block one by one, from the bottom up, as pivots and at each stage we will perform elementary column operations to eliminate non-zero entries to the left of the pivot. By our remark on upper triangularness, these operations only affect the first $|R|-1$ columns. After elimination, the upper left $(|R|-1)\times(|R|-1)$ block $B$ is such that its determinant differs from $\det M_{r\to e}$ only by a monomial factor, namely the product of the pivots. (In fact the factor is the hypertree in $(E,R)$, written as a monomial, which is associated to the chosen arborescence/Tutte matching.) Therefore it suffices to show that for each $j=1,\ldots,|R|-1$, the $j$'th column of $B$ (indexed by the triangle $t_j$) coincides with $U_j^{-1}$ times the $j$'th column (indexed by the word $W_j$) of the Turaev torsion.

\emph{Why the two determinants are equal:} We will start with a rough description of the effect of our elimination steps and fill in some details later. It suffices to concentrate on just one arbitrary column, say the $j$'th. In the original enhanced adjacency matrix, the $j$'th column contains three non-zero entries: $1$ in the $j$'th row which corresponds to the summand $U_j$ in $\partial W_j/\partial r_j$; another $1$ in the row indexed by $v^*$; and $r_j$ in the row indexed by $e^*$. If $e^*$ or $v^*$ is a root, then one of the latter two is missing. During the elimination process, the number of non-zero entries in the bottom $n-|R|$ positions never increases: it stays two for a while, then goes down to one and eventually to zero. (Or else, it stays one for a while and then becomes zero.)

It is best to think of what is happening in the $j$'th column as two superimposed, \emph{left and right processes} $LP$ and $RP$, one for each section of $\lambda_j$ before and after the edge $v^*e^*$. (If $e^*$ or $v^*$ is a root, then the corresponding process is empty.) Each process has a non-zero monomial entry, called the \emph{can}, in one of the bottom $n-|R|$ positions which it moves gradually upward until it disappears (as in `kicking the can down the road'). The can is not constant in the process. In each elimination step, if it affects the $j$'th column, then one of the two processes takes a step. In that step, the presence of the two (or one) other non-zero entries in the column of the pivot results in the following:
\begin{enumerate}
\item\label{can} A new non-zero monomial entry (can) is created in one of the lower $n-|R|$ positions of the $j$'th column. It is higher up than the entry which is eliminated in the step. If the column of the pivot belongs to a white triangle which is adjacent to $e_0$ or $v_0$, then this development does not occur and the process terminates.
\item\label{contribution} A new summand (contribution) is added to one of the entries in the upper $|R|-1$ positions of the $j$'th row. If the column of the pivot belongs to a white triangle adjacent to $r_0$, then no contribution occurs.
\end{enumerate}
It is also not hard to see that the sequence of pivots for $LP$ and $RP$ starts in the rows indexed by $v^*$ and $e^*$, respectively, and that the rows of other pivots are indexed by the emerald and violet points along $\lambda_j$ that occur as we move toward $v_0$ along the respective paths. In particular, emerald pivots (equal to some $r_i$) and violet pivots (equal to $1$) alternate in $LP$ as well as in $RP$. If at either end of $\lambda_j$ the point adjacent to $v_0$ is $e_0$, then the corresponding process ends after the pivot which is in the row of the other violet neighbor of $e_0$ along $\lambda_j$. It is also possible for the two processes to merge after a while as edges and white triangles at the beginning and end of $\lambda_j$ coincide. After this happens, the cans of $LP$ and $RP$ occur in the same position as a sum of two monomials.

Finally, we need to examine the monomials (cans and contributions) in \eqref{can} and \eqref{contribution} above in more concrete terms. Please refer to Figure \ref{fig:snake} for notation. Let us first consider the left process $LP$ and the stage when the pivot is the $(k,k)$-entry $y$, eliminating the can $C^y$ in the position $(k,j)$. The new can is $-y^{-1}C^y$ and it is in the row indexed by $p_{k'}$. A contribution of $-y^{-1}C^y$ is added to the $j$'th entry in the row indexed by $y$. Let now the pivot be the $(l,l)$-entry $1$ and denote the corresponding can by $C^x$. In this case the new can is $-xC^x$ and it is in the row indexed by $p_{l'}$. A contribution of $-C^x$ is added to the $j$'th entry in the row indexed by $x$.

From this it is clear that the sequence of cans for $LP$ is $1$, $-a$, $b^{-1}a$, and so on so that $C^y$ is $-1$ times the product (with alternating exponents) of red points (labels) along $\lambda_j$ between $p_k$ and $v^*$ and $C^x$ is described similarly but without the negative sign. Then if we look at the contributions $-y^{-1}C^y$ and $-C^x$ of the previous paragraph, we see that if we multiply them by $U_j$, we get $W^y$ and $-W^xx^{-1}$, respectively, just as we expected and exactly where we expected them.

A similar analysis applies to the right process $RP$. There, the first can is $r_j$ in the row of $e^*$, the next one is $-c^{-1}r_j$ (note that the first pivot is $c$) in the row indexed by $w$ and so on. The entry $1$ in the $(j,j)$-position plays the role of first contribution, then the second one is $-c^{-1}r_j$ in the row indexed by $c$ and so on. After multiplying with $U_j$, these also conform to our expected values.
\end{proof}

Having established Theorem \ref{thm:main}, Corollary \ref{cor:samesize} is now immediate: $|S(F_G)|$ is the number of hypertrees in $(E,R)$, but that is the same (by Proposition \ref{pro:det}) as $|\det M|$ of the `un-enhanced' adjacency matrix $M$ of the trinity. That number (the number of Tutte matchings) is clearly color-independent.

Let us recall that if $F$ is of minimal genus, then $\sum_{\mathfrak s}\chi(SFH(M_F,L,\mathfrak s))$ is the leading coefficient of the Alexander polynomial $\Delta_L$, cf.\ \cite[Lemma 6.4]{decat} and \cite[Theorem 1.5]{sfhdecomp}. (Note that $\Delta$ and Conway's version $\nabla$ have the same leading coefficient.) Thus, in our class of examples, $|S(F_G)|$ is the leading coefficient of $\Delta_{L_G},$ and we also obtain the following result.

\begin{corollary}\label{cor:alex}
The determinant of the adjacency matrix, i.e., the number of Tutte matchings, is the leading coefficient in the Alexander polynomial of any of the three alternating links associated to the trinity.
\end{corollary}

This last fact has a more direct proof \cite{km} using the Tree Trinity Theorem \cite{tutte} and Kauffman's state model \cite{kauff} for the Alexander polynomial. We now prove Corollary \ref{Cor:convex}.

\pagebreak

\begin{proof}
First, suppose \(F_G\) is a standard Seifert surface of a special
alternating link \(L_G\), so that \(S(F_G) \cong Q_{(E,R)} \).
By \cite[Theorem 3.4]{hipertutte}, \(Q_{(E,R)}\) is cut out by a
system of linear inequalities in \(\Z^R\), so it is clearly convex.

Next, suppose \(F\) is a standard Seifert surface of a homogeneous link. It
is well known that any such \(F\)
is a Murasugi sum of standard Seifert surfaces
for special alternating links. Now \eqref{eqn:Murasugi}, along with the fact that
the product of convex sets is convex, proves the desired result.
\end{proof}



As another application of Theorem \ref{thm:main}, we analyze the following situation. Let the edge $\varepsilon$ of the $2$-connected plane bipartite graph $G$ be adjacent to the regions $r_1$ and $r_2$. Let $G'$ be the graph in which $\varepsilon$ is replaced with a path of three edges. Then the handlebodies $M_G$ and $M_{G'}$ have a natural identification so that we can think of the supports $S(F_G)$ and $S(F_{G'})$ as subsets of the same affine space over $H_1(M_G;\Z)$. In that group, let $\varepsilon^*$ denote the homology class of a meridian (in $S^3$) of $\varepsilon$.

\begin{proposition}\label{pro:lengthen}
With the above conventions, $S(F_{G'})\cong S(F_G)\cup(S(F_G)+\varepsilon^*)$, meaning that the two sides differ by a translation.
\end{proposition}

For the convex hulls of the supports, the statement is that $\conv(S(F_{G'}))$ is the Minkowski sum of $\conv(S(F_G))$ and a certain line segment. Note that our main example is an instance of $G'$, which explains the `elongated' shape of the polytope in Figure \ref{fig:polytope} in the direction from $r_1^3$ to $r_2^3$.

\begin{proof}
We will prove the relevant claim on sets of hypertrees and leave it to the reader to work out the various identifications. Let $E$ and $V$ denote the color classes of $G,$ and let $E'=E\cup\{e^*\}$ and $V'=V\cup\{v^*\}$ be the corresponding color classes in $G'$. Both $G$ and $G'$ have the same set $R$ of regions. Then, we claim that
\[Q_{(E',R)}=(Q_{(E,R)}+\mathbf i_{\{r_1\}})\cup(Q_{(E,R)}+\mathbf i_{\{r_2\}}),\]
where the $\mathbf i_{\{r_i\}}$ are standard generators in $\Z^R$ (indicator functions of singleton sets).

To see this, note that any hypertree $\mathbf f$ in $(E,R)$ induces two hypertrees in $(E',R)$ via increasing by $1$ the value of $\mathbf f$ either at $r_1$ or at $r_2$. Indeed, these have obvious spanning trees (in $\bip(E',R)$) realizing them which are built from a realization (in $\bip(E,R)$) of $\mathbf f$ by adding either the edge $r_1e^*$ or the edge $r_2e^*$.

To establish the converse, it suffices to show that any hypertree $\mathbf g$ in $(E',R)$ has a realization that is of valence one at $e^*$. If the realization $\Gamma$ of $\mathbf g$ is not such, then $e^*$ is a valence two point in it. Let $e$ denote the emerald endpoint of $\varepsilon$. Then $r_1e$ and $r_2e$ are edges in $\bip(E',R),$ but at most one of them can be an edge in $\Gamma$. If one, say $r_1e$, is in $\Gamma$, then remove $r_2e^*$ from $\Gamma$ and replace it with $r_2e$ to get the realization of $\mathbf g$ with the desired property. If neither $r_1e$ nor $r_2e$ is in $\Gamma$, then add $r_1e$ to $\Gamma$ and kill the resulting cycle by removing its other edge adjacent to $r_1$. If this was $r_1e^*$, we are done; otherwise, apply the previous step.
\end{proof}


\pagebreak

\section{A sample calculation}\label{sec:ex}

\begin{figure}[htbp]
\labellist
\footnotesize
\pinlabel $r_0$ at 366 98
\pinlabel $r_1$ at 259 422
\pinlabel $r_2$ at 601 368
\pinlabel $r_3$ at 528 633
\pinlabel $v_0$ at 276 602
\pinlabel $e_0$ at 834 529
\tiny
\pinlabel V at 387 462
\pinlabel II at 351 314
\pinlabel I at 689 206
\pinlabel VI at 621 548
\pinlabel VII at 419 530
\pinlabel IV at 365 391
\pinlabel III at 347 242
\endlabellist
   \centering
   \includegraphics[width=\linewidth]{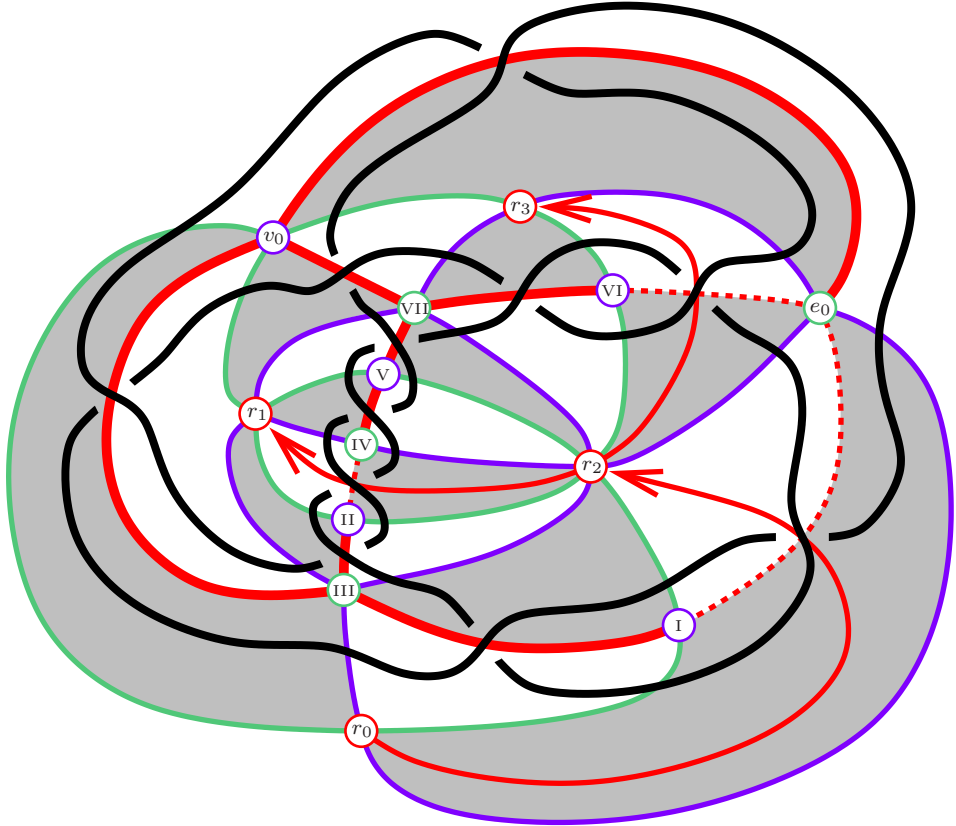}
   \caption{The spanning arborescence $A\subset G_R^*$, its dual tree $\Gamma\subset G_R$, and a deconstruction order for the latter.}
   \label{fig:compute}
\end{figure}

We will show explicit reduction steps on the enhanced adjacency matrix in our main running example. The computation will use the spanning arborescence (of $G_R^*$) indicated in Figure \ref{fig:compute}. We also fix the deconstruction order $v_3<v_2<e_3<e_2<v_1<v_4<e_1$ as shown in the figure. For better visibility, we have rearranged the rows and columns of the matrix \eqref{adjmatrix} to reflect this order \emph{and} the Tutte matching determined by the arborescence: it now appears as the main diagonal. (Thus we see that the chosen arborescence corresponds to the hypertree/monomial $r_1r_2r_3$ of Example~\ref{ex:det}.) We briefly indicate the $p$-labels that appear in the proof of Theorem~\ref{thm:main}, but we \emph{did not} re-label the white triangles. The signs of all expansion terms in the determinant are still positive, but this is not just by coincidence any more, rather by design. Note how the bottom right $7\times7$ submatrix is already upper triangular at the start.

In each of the seven steps below, we use the encircled pivot to clear out non-zero entries to its left. Each column of the top left $3\times3$ block that remains at the end is equal to a monomial times the corresponding column of the determinant \eqref{turaev}. As the computation in Example \ref{ex:paths} was carried out using the procedure in the first half of the proof of Theorem \ref{thm:main}, this indeed illustrates the correctness of our method.

\[\bordermatrix
{~&\text{\footnotesize $t_4$}&\text{\footnotesize $t_{10}$}&\text{\footnotesize $t_9$}&\text{\footnotesize $t_1$}&\text{\footnotesize $t_8$}&\text{\footnotesize $t_2$}&\text{\footnotesize $t_7$}&\text{\footnotesize $t_3$}&\text{\footnotesize $t_6$}&\text{\footnotesize $t_5$}\cr
\makebox[25pt][r]{\text{\footnotesize $r_1$}}&1&0&0&0&0&1&0&1&0&0\cr
\makebox[25pt][r]{\text{\footnotesize $r_2$}}&0&1&0&0&1&0&1&0&1&0\cr
\makebox[25pt][r]{\text{\footnotesize $r_3$}}&0&0&1&1&0&0&0&0&0&0\cr
\makebox[25pt][r]{\text{\footnotesize $p_4=e_1$}}&0&0&0&r_3&r_2&r_1&0&0&0&0\cr
\makebox[25pt][r]{\text{\footnotesize $p_5=v_4$}}&0&0&1&0&1&0&0&0&0&0\cr
\makebox[25pt][r]{\text{\footnotesize $p_6=v_1$}}&0&0&0&0&0&1&1&0&0&0\cr
\makebox[25pt][r]{\text{\footnotesize $p_7=e_2$}}&r_1&0&0&0&0&0&r_2&0&0&0\cr
\makebox[25pt][r]{\text{\footnotesize $p_8=e_3$}}&0&0&0&0&0&0&0&r_1&r_2&r_0\cr
\makebox[25pt][r]{\text{\footnotesize $p_9=v_2$}}&1&0&0&0&0&0&0&0&1&0\cr
\makebox[25pt][r]{\text{\footnotesize $p_{10}=v_3$}}&0&1&0&0&0&0&0&0&0&\begin{picture}(0,0)\put(0,0){\makebox(0,0)[b]{$1$}}\put(0,3){\circle{13}}\end{picture}\cr};\]

\[\bordermatrix
{~&\text{\footnotesize $t_4$}&\text{\footnotesize $t_{10}$}&\text{\footnotesize $t_9$}&\text{\footnotesize $t_1$}&\text{\footnotesize $t_8$}&\text{\footnotesize $t_2$}&\text{\footnotesize $t_7$}&\text{\footnotesize $t_3$}&\text{\footnotesize $t_6$}&\text{\footnotesize $t_5$}\cr
\text{\footnotesize $r_1$}&1&0&0&0&0&1&0&1&0&0\cr
\text{\footnotesize $r_2$}&0&1&0&0&1&0&1&0&1&0\cr
\text{\footnotesize $r_3$}&0&0&1&1&0&0&0&0&0&0\cr
\text{\footnotesize $e_1$}&0&0&0&r_3&r_2&r_1&0&0&0&0\cr
\text{\footnotesize $v_4$}&0&0&1&0&1&0&0&0&0&0\cr
\text{\footnotesize $v_1$}&0&0&0&0&0&1&1&0&0&0\cr
\text{\footnotesize $e_2$}&r_1&0&0&0&0&0&r_2&0&0&0\cr
\text{\footnotesize $e_3$}&0&-r_0&0&0&0&0&0&r_1&r_2&r_0\cr
\text{\footnotesize $v_2$}&1&0&0&0&0&0&0&0&\begin{picture}(0,0)\put(0,0){\makebox(0,0)[b]{$1$}}\put(0,3){\circle{13}}\end{picture}&0\cr
\text{\footnotesize $v_3$}&0&0&0&0&0&0&0&0&0&\begin{picture}(0,0)\put(0,0){\makebox(0,0)[b]{$1$}}\color{gray}\put(0,3){\circle{13}}\end{picture}\cr};\]

\[\bordermatrix
{~&\text{\footnotesize $t_4$}&\text{\footnotesize $t_{10}$}&\text{\footnotesize $t_9$}&\text{\footnotesize $t_1$}&\text{\footnotesize $t_8$}&\text{\footnotesize $t_2$}&\text{\footnotesize $t_7$}&\text{\footnotesize $t_3$}&\text{\footnotesize $t_6$}&\text{\footnotesize $t_5$}\cr
\text{\footnotesize $r_1$}&1&0&0&0&0&1&0&1&0&0\cr
\text{\footnotesize $r_2$}&-1&1&0&0&1&0&1&0&1&0\cr
\text{\footnotesize $r_3$}&0&0&1&1&0&0&0&0&0&0\cr
\text{\footnotesize $e_1$}&0&0&0&r_3&r_2&r_1&0&0&0&0\cr
\text{\footnotesize $v_4$}&0&0&1&0&1&0&0&0&0&0\cr
\text{\footnotesize $v_1$}&0&0&0&0&0&1&1&0&0&0\cr
\text{\footnotesize $e_2$}&r_1&0&0&0&0&0&r_2&0&0&0\cr
\text{\footnotesize $e_3$}&-r_2&-r_0&0&0&0&0&0&\begin{picture}(0,0)\put(0,0){\makebox(0,0)[b]{$r_1$}}\put(0,3){\circle{13}}\end{picture}&r_2&r_0\cr
\text{\footnotesize $v_2$}&0&0&0&0&0&0&0&0&\begin{picture}(0,0)\put(0,0){\makebox(0,0)[b]{$1$}}\color{gray}\put(0,3){\circle{13}}\end{picture}&0\cr
\text{\footnotesize $v_3$}&0&0&0&0&0&0&0&0&0&\begin{picture}(0,0)\put(0,0){\makebox(0,0)[b]{$1$}}\color{gray}\put(0,3){\circle{13}}\end{picture}\cr};\]

\[\bordermatrix
{~&\text{\footnotesize $t_4$}&\text{\footnotesize $t_{10}$}&\text{\footnotesize $t_9$}&\text{\footnotesize $t_1$}&\text{\footnotesize $t_8$}&\text{\footnotesize $t_2$}&\text{\footnotesize $t_7$}&\text{\footnotesize $t_3$}&\text{\footnotesize $t_6$}&\text{\footnotesize $t_5$}\cr
\text{\footnotesize $r_1$}&1+r_1^{-1}r_2&r_1^{-1}r_0&0&0&0&1&0&1&0&0\cr
\text{\footnotesize $r_2$}&-1&1&0&0&1&0&1&0&1&0\cr
\text{\footnotesize $r_3$}&0&0&1&1&0&0&0&0&0&0\cr
\text{\footnotesize $e_1$}&0&0&0&r_3&r_2&r_1&0&0&0&0\cr
\text{\footnotesize $v_4$}&0&0&1&0&1&0&0&0&0&0\cr
\text{\footnotesize $v_1$}&0&0&0&0&0&1&1&0&0&0\cr
\text{\footnotesize $e_2$}&r_1&0&0&0&0&0&\begin{picture}(0,0)\put(0,0){\makebox(0,0)[b]{$r_2$}}\put(0,3){\circle{13}}\end{picture}&0&0&0\cr
\text{\footnotesize $e_3$}&0&0&0&0&0&0&0&\begin{picture}(0,0)\put(0,0){\makebox(0,0)[b]{$r_1$}}\color{gray}\put(0,3){\circle{13}}\end{picture}&r_2&r_0\cr
\text{\footnotesize $v_2$}&0&0&0&0&0&0&0&0&\begin{picture}(0,0)\put(0,0){\makebox(0,0)[b]{$1$}}\color{gray}\put(0,3){\circle{13}}\end{picture}&0\cr
\text{\footnotesize $v_3$}&0&0&0&0&0&0&0&0&0&\begin{picture}(0,0)\put(0,0){\makebox(0,0)[b]{$1$}}\color{gray}\put(0,3){\circle{13}}\end{picture}\cr};\]

\[\bordermatrix
{~&\text{\footnotesize $t_4$}&\text{\footnotesize $t_{10}$}&\text{\footnotesize $t_9$}&\text{\footnotesize $t_1$}&\text{\footnotesize $t_8$}&\text{\footnotesize $t_2$}&\text{\footnotesize $t_7$}&\text{\footnotesize $t_3$}&\text{\footnotesize $t_6$}&\text{\footnotesize $t_5$}\cr
\text{\footnotesize $r_1$}&1+r_1^{-1}r_2&r_1^{-1}r_0&0&0&0&1&0&1&0&0\cr
\text{\footnotesize $r_2$}&-1-r_2^{-1}r_1&1&0&0&1&0&1&0&1&0\cr
\text{\footnotesize $r_3$}&0&0&1&1&0&0&0&0&0&0\cr
\text{\footnotesize $e_1$}&0&0&0&r_3&r_2&r_1&0&0&0&0\cr
\text{\footnotesize $v_4$}&0&0&1&0&1&0&0&0&0&0\cr
\text{\footnotesize $v_1$}&-r_2^{-1}r_1&0&0&0&0&\begin{picture}(0,0)\put(0,0){\makebox(0,0)[b]{$1$}}\put(0,3){\circle{13}}\end{picture}&1&0&0&0\cr
\text{\footnotesize $e_2$}&0&0&0&0&0&0&\begin{picture}(0,0)\put(0,0){\makebox(0,0)[b]{$r_2$}}\color{gray}\put(0,3){\circle{13}}\end{picture}&0&0&0\cr
\text{\footnotesize $e_3$}&0&0&0&0&0&0&0&\begin{picture}(0,0)\put(0,0){\makebox(0,0)[b]{$r_1$}}\color{gray}\put(0,3){\circle{13}}\end{picture}&r_2&r_0\cr
\text{\footnotesize $v_2$}&0&0&0&0&0&0&0&0&\begin{picture}(0,0)\put(0,0){\makebox(0,0)[b]{$1$}}\color{gray}\put(0,3){\circle{13}}\end{picture}&0\cr
\text{\footnotesize $v_3$}&0&0&0&0&0&0&0&0&0&\begin{picture}(0,0)\put(0,0){\makebox(0,0)[b]{$1$}}\color{gray}\put(0,3){\circle{13}}\end{picture}\cr};\]

\[\bordermatrix
{~&\text{\footnotesize $t_4$}&\text{\footnotesize $t_{10}$}&\text{\footnotesize $t_9$}&\text{\footnotesize $t_1$}&\text{\footnotesize $t_8$}&\text{\footnotesize $t_2$}&\text{\footnotesize $t_7$}&\text{\footnotesize $t_3$}&\text{\footnotesize $t_6$}&\text{\footnotesize $t_5$}\cr
\text{\footnotesize $r_1$}&1+r_1^{-1}r_2+r_2^{-1}r_1&r_1^{-1}r_0&0&0&0&1&0&1&0&0\cr
\text{\footnotesize $r_2$}&-1-r_2^{-1}r_1&1&0&0&1&0&1&0&1&0\cr
\text{\footnotesize $r_3$}&0&0&1&1&0&0&0&0&0&0\cr
\text{\footnotesize $e_1$}&r_2^{-1}r_1^2&0&0&r_3&r_2&r_1&0&0&0&0\cr
\text{\footnotesize $v_4$}&0&0&1&0&\begin{picture}(0,0)\put(0,0){\makebox(0,0)[b]{$1$}}\put(0,3){\circle{13}}\end{picture}&0&0&0&0&0\cr
\text{\footnotesize $v_1$}&0&0&0&0&0&\begin{picture}(0,0)\put(0,0){\makebox(0,0)[b]{$1$}}\color{gray}\put(0,3){\circle{13}}\end{picture}&1&0&0&0\cr
\text{\footnotesize $e_2$}&0&0&0&0&0&0&\begin{picture}(0,0)\put(0,0){\makebox(0,0)[b]{$r_2$}}\color{gray}\put(0,3){\circle{13}}\end{picture}&0&0&0\cr
\text{\footnotesize $e_3$}&0&0&0&0&0&0&0&\begin{picture}(0,0)\put(0,0){\makebox(0,0)[b]{$r_1$}}\color{gray}\put(0,3){\circle{13}}\end{picture}&r_2&r_0\cr
\text{\footnotesize $v_2$}&0&0&0&0&0&0&0&0&\begin{picture}(0,0)\put(0,0){\makebox(0,0)[b]{$1$}}\color{gray}\put(0,3){\circle{13}}\end{picture}&0\cr
\text{\footnotesize $v_3$}&0&0&0&0&0&0&0&0&0&\begin{picture}(0,0)\put(0,0){\makebox(0,0)[b]{$1$}}\color{gray}\put(0,3){\circle{13}}\end{picture}\cr};\]

\[\bordermatrix
{~&\text{\footnotesize $t_4$}&\text{\footnotesize $t_{10}$}&\text{\footnotesize $t_9$}&\text{\footnotesize $t_1$}&\text{\footnotesize $t_8$}&\text{\footnotesize $t_2$}&\text{\footnotesize $t_7$}&\text{\footnotesize $t_3$}&\text{\footnotesize $t_6$}&\text{\footnotesize $t_5$}\cr
\text{\footnotesize $r_1$}&1+r_1^{-1}r_2+r_2^{-1}r_1&r_1^{-1}r_0&0&0&0&1&0&1&0&0\cr
\text{\footnotesize $r_2$}&-1-r_2^{-1}r_1&1&-1&0&1&0&1&0&1&0\cr
\text{\footnotesize $r_3$}&0&0&1&1&0&0&0&0&0&0\cr
\text{\footnotesize $e_1$}&r_2^{-1}r_1^2&0&-r_2&\begin{picture}(0,0)\put(0,0){\makebox(0,0)[b]{$r_3$}}\put(0,3){\circle{13}}\end{picture}&r_2&r_1&0&0&0&0\cr
\text{\footnotesize $v_4$}&0&0&0&0&\begin{picture}(0,0)\put(0,0){\makebox(0,0)[b]{$1$}}\color{gray}\put(0,3){\circle{13}}\end{picture}&0&0&0&0&0\cr
\text{\footnotesize $v_1$}&0&0&0&0&0&\begin{picture}(0,0)\put(0,0){\makebox(0,0)[b]{$1$}}\color{gray}\put(0,3){\circle{13}}\end{picture}&1&0&0&0\cr
\text{\footnotesize $e_2$}&0&0&0&0&0&0&\begin{picture}(0,0)\put(0,0){\makebox(0,0)[b]{$r_2$}}\color{gray}\put(0,3){\circle{13}}\end{picture}&0&0&0\cr
\text{\footnotesize $e_3$}&0&0&0&0&0&0&0&\begin{picture}(0,0)\put(0,0){\makebox(0,0)[b]{$r_1$}}\color{gray}\put(0,3){\circle{13}}\end{picture}&r_2&r_0\cr
\text{\footnotesize $v_2$}&0&0&0&0&0&0&0&0&\begin{picture}(0,0)\put(0,0){\makebox(0,0)[b]{$1$}}\color{gray}\put(0,3){\circle{13}}\end{picture}&0\cr
\text{\footnotesize $v_3$}&0&0&0&0&0&0&0&0&0&\begin{picture}(0,0)\put(0,0){\makebox(0,0)[b]{$1$}}\color{gray}\put(0,3){\circle{13}}\end{picture}\cr};\]

\[\bordermatrix
{~&\text{\footnotesize $t_4$}&\text{\footnotesize $t_{10}$}&\text{\footnotesize $t_9$}&\text{\footnotesize $t_1$}&\text{\footnotesize $t_8$}&\text{\footnotesize $t_2$}&\text{\footnotesize $t_7$}&\text{\footnotesize $t_3$}&\text{\footnotesize $t_6$}&\text{\footnotesize $t_5$}\cr
\text{\footnotesize $r_1$}&1+r_1^{-1}r_2+r_2^{-1}r_1&r_1^{-1}r_0&0&0&0&1&0&1&0&0\cr
\text{\footnotesize $r_2$}&-1-r_2^{-1}r_1&1&-1&0&1&0&1&0&1&0\cr
\text{\footnotesize $r_3$}&-r_3^{-1}r_2^{-1}r_1^2&0&1+r_3^{-1}r_2&1&0&0&0&0&0&0\cr
\text{\footnotesize $e_1$}&0&0&0&\begin{picture}(0,0)\put(0,0){\makebox(0,0)[b]{$r_3$}}\color{gray}\put(0,3){\circle{13}}\end{picture}&r_2&r_1&0&0&0&0\cr
\text{\footnotesize $v_4$}&0&0&0&0&\begin{picture}(0,0)\put(0,0){\makebox(0,0)[b]{$1$}}\color{gray}\put(0,3){\circle{13}}\end{picture}&0&0&0&0&0\cr
\text{\footnotesize $v_1$}&0&0&0&0&0&\begin{picture}(0,0)\put(0,0){\makebox(0,0)[b]{$1$}}\color{gray}\put(0,3){\circle{13}}\end{picture}&1&0&0&0\cr
\text{\footnotesize $e_2$}&0&0&0&0&0&0&\begin{picture}(0,0)\put(0,0){\makebox(0,0)[b]{$r_2$}}\color{gray}\put(0,3){\circle{13}}\end{picture}&0&0&0\cr
\text{\footnotesize $e_3$}&0&0&0&0&0&0&0&\begin{picture}(0,0)\put(0,0){\makebox(0,0)[b]{$r_1$}}\color{gray}\put(0,3){\circle{13}}\end{picture}&r_2&r_0\cr
\text{\footnotesize $v_2$}&0&0&0&0&0&0&0&0&\begin{picture}(0,0)\put(0,0){\makebox(0,0)[b]{$1$}}\color{gray}\put(0,3){\circle{13}}\end{picture}&0\cr
\text{\footnotesize $v_3$}&0&0&0&0&0&0&0&0&0&\begin{picture}(0,0)\put(0,0){\makebox(0,0)[b]{$1$}}\color{gray}\put(0,3){\circle{13}}\end{picture}\cr}.\]


\begin{thebibliography}{00}

\bibitem{irida}Irida Altman, \textit{Sutured Floer homology distinguishes between Seifert surfaces}, arXiv:1012.5904.

\bibitem{banks}Jessica E. Banks, \textit{Minimal Genus Seifert Surfaces for Alternating Links}, arXiv:1106.3180.

\bibitem{berman}
K.\ A.\ Berman, \textit{A proof of Tutte's trinity theorem and a new determinant formula}, SIAM J.\ Algebraic Discrete Methods 1 (1980), no.\ 1, 64--69.

\bibitem{cromwell}
P. R. Cromwell, \textit{Homogeneous links}, J. London Math.\ Soc.\ 39 (1989), no.\ 3, 535--552.

\bibitem{decat}S.\ Friedl, A.\ Juh\'asz, and J.\ Rasmussen,
  \textit{The decategorification of sutured Floer homology}, arXiv:0903.5287.

\bibitem{gabai} D.\ Gabai, \textit{Foliations and the topology of 3-manifolds},
  J.\ Diff.\ Geom.\ 18 (1983), 445--503.

\bibitem{hisa}M. Hirasawa and M. Sakuma, \textit{Minimal genus Seifert surfaces for alternating links}, in KNOTS '96 (Tokyo), pp.\ 383--394, World Sci.\ Publ., River Edge, NJ, 1997.

\bibitem{jaeger}F.\ Jaeger, \textit{Tutte polynomials and link polynomials}, Proc.\ Amer.\ Math.\ Soc.\ 103 (1988), no.\ 2, 647--654.

\bibitem{sfh}A.\ Juh\'asz, \textit{Holomorphic discs and sutured
    manifolds}, Algebr.\ Geom.\ Topol.\ 6 (2006), 1429--1457.

\bibitem{sfhdecomp} A.\ Juh\'asz, \textit{Floer homology and surface
  decompositions}, Geom.\ Topol.\ 12 (2008), 299--350.

\bibitem{sfhseifert} A.\ Juh\'asz, \textit{Knot Floer homology and
  Seifert surfaces}, Algebr.\ Geom.\ Topol.\ 8 (2008), 603--608.

\bibitem{sfhpolytope}A. Juh\'asz, \textit{The sutured Floer homology polytope}, Geom.\ Topol.\ 14 (2010), no.\ 3, 1303--1354.

\bibitem{hipertutte}T.\ K\'alm\'an, \textit{A version of Tutte's polynomial for hypergraphs}, arXiv:1103.1057.

\bibitem{km}T.\ K\'alm\'an and H.\ Murakami, \textit{Root polytopes and the Homfly polynomial of alternating links}, in progress.


\bibitem{kauff}L. H. Kauffman, \textit{Formal Knot Theory}, Mathematical Notes 30, Princeton Univ.\ Press, Princeton, NJ, 1983.

\bibitem{mprz}K.\ Murasugi and J.\ Przytycki, \textit{The skein polynomial of a planar star product of two links}, Math.\ Proc.\ Camb.\ Phil.\ Soc.\ 106 (1989), 273--276.

\bibitem{post}A.\ Postnikov, \textit{Permutohedra, Associahedra, and Beyond}, Int.\ Math.\ Res.\ Not.\ 2009, no.\ 6, 1026--1106.


\bibitem{tutte}W.\ T.\ Tutte, \textit{The dissection of equilateral triangles into equilateral triangles}, Proc.\ Cambridge Philos.\ Soc.\ 44 (1948), 463--482.

\bibitem{Webb} R. Webb, \textit{Alexander polynomials of sutured handlebodies}, Trinity College summer research report, 2009.

\end{thebibliography}
\end{document}